\documentclass[letterpaper,10pt]{article} 

\usepackage{osameet3} 

\usepackage{amsmath,amssymb}
\usepackage{subfigure}
\usepackage{wrapfig}
\usepackage{algpseudocode,algorithm,algorithmicx}
\algrenewcommand\algorithmicrequire{\textbf{Initialize:}}

\begin{document}

\title{Multi-dimensional Joint Probabilistic and Geometric Shaping Strategies for Nonlinear Fiber Communications}

\author{Metodi P. Yankov}
\address{Department of Photonics Engineering, Tecgnical University of Denmark, 2800 Kgs. Lyngby, Denmark}
\email{meya@fotonik.dtu.dk}

\copyrightyear{2020}

\begin{abstract}
Constellation shaping is reviewed and revised for a WDM unrepeated system with high spectral efficiency. It is shown that for a constellation size-constrained system, previous optimization techniques can be highly sub-optimal, and a superior optimization algorithm is proposed.
\end{abstract}

\ocis{060.1660, 060.2330.}

\section{Introduction}
Constellation shaping is a requirement for digital communication systems to achieve the channel capacity \cite{ThomasCover}. Shaping inevitably entails optimization of the transmission alphabet to match the channel conditions. For discrete, memoryless channels (DMCs), optimization strategies producing optimal solutions have been known since the 70s \cite{ThomasCover}, with their optimal modifications for e.g. the additive white Gaussian noise (AWGN) channel \cite{Varnica}. Theoretically, joint geometric (on the constellation points positions) and probabilistic (on their probability of occurrence, collectively referred to as probability mass function (PMF)) shaping is required to achieve the channel capacity. In the above-mentioned DMC cases and in the case of highly-spectral efficient modulation formats, either is typically sufficient. 
Recently, due to the ever-increasing data rate demand, optical communication systems have been pushed to operate ever-closer to channel capacity, inevitably facing the problem of shaping. Several practical solutions have already been proposed, e.g. \cite{Buchali}, adopting different variants of the probabilistic amplitude shaping scheme for AWGN channels \cite{Bocherer}. Typical optical receivers operate under the memoryless Gaussian channel assumption, for which optimization with the above-mentioned schemes produces reasonable shaping gains. Increasing the data rate in such systems entails entering the nonlinear region of transmission, where the above-mentioned Gaussian channel assumption becomes less and less accurate. In particular, through the nonlinearities, the channel becomes dependent on the input PMF. Furthermore, in such cases, the near-optimality assumption of independent probabilistic and geometric shaping should be revisited. Hence, new constellation optimization strategies are needed to unlock the nonlinear region of transmission. 

This paper proposes such a solution, targeting the nonlinear region of transmission. Optimization is limited to 4D, i.e. dual-polarization transmission, however, the optimization method is general to multi-dimensional (MD) signals, e.g. joint sub-carrier constellation shaping or joint time-slot constellation optimization.

\section{Optimization algorithm}
\begin{wrapfigure}{r}{0.6\textwidth}
  \begin{center}
  \vspace{-0.7cm}
    \includegraphics[width=9.5cm]{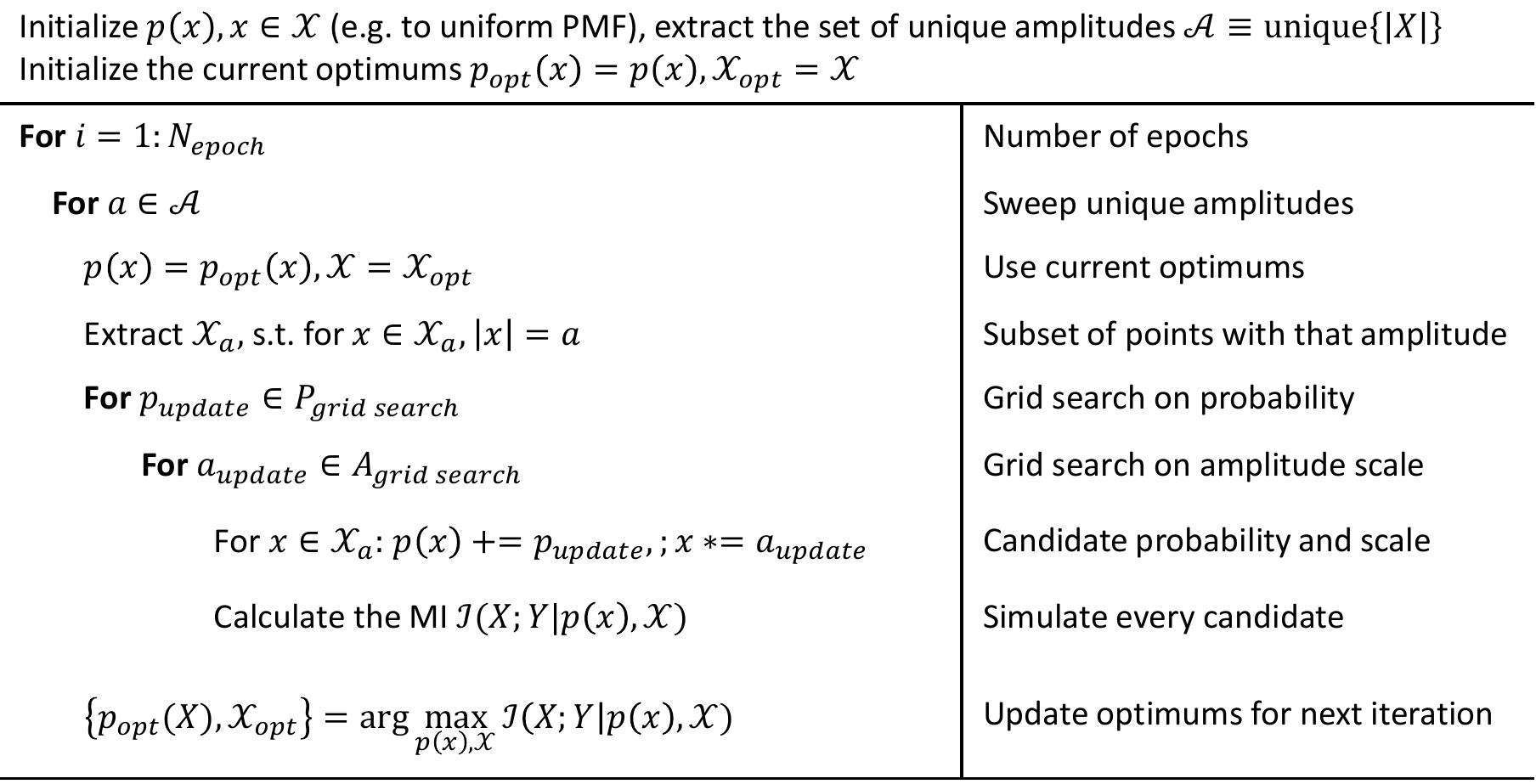}
  \end{center}
  \vspace{-0.5cm}
\end{wrapfigure}

In order to optimize the PMF of the constellation, a greedy approach is taken. Greedy algorithms reduce the complexity of high-dimensional optimization problems by performing local optimization on a subset of parameters and iterating between local optimizations a number of times (referred to as epochs). In the proposed algorithm, local optimization is performed on the probabilities and geometric scale of points, that share the same MD amplitude. I.e., joint geometric and probabilistic amplitude optimization is performed on a base 4D constellation ${\cal{X}}$, e.g. ${\cal{X}}=64^2$QAM. The set of unique amplitudes ${\cal{A}}$ is then swept, and epochs are performed until convergence. The local optimization is based on a 2D grid search - probability dimension and scale dimension - around the current local optimum. A summary of the algorithm is given above. It is noted that the grid search contains the zero-probability points, i.e., points are allowed to be pruned away. We observed convergence of the algorithm after $\approx$ 5 epochs. The algorithm is agnostic to the channel under optimization. For example, for an AWGN channel, the algorithm successfuly converges to the known-to-be-near-optimal Maxwell-Boltzmann (MB) PMF. The only requirement is that the channel can be simulated. 

\subsection{System under test}
The simulation setup is given in Fig.~1(a). An unrepeated system for which the optimum launch power is relatively high is used to exemplify the algorithm. In order to emulate such a system, a 5x30 GBaud, 50 GHz spacing wavelength division multiplexed (WDM) signal is generated, the split-step Fourier method is used to simulate 250 km of standard, single mode fiber propagation, and an erbium doped fiber amplifier (EDFA) (noise figure of 5 dB) is used at the receiver that sets the received power to 0 dBm. The chromatic dispersion is then compensated in the frequency domain, matched filtering is performed to extract the central channel, and the mismatched decoding technique \cite{Yankov} with a 4D Gaussian auxiliary channel \cite{Fehenberger} is used to estimate the mutual information (MI), representing an achievable information rate (AIR).


\subsection{Reference optimization strategies}
The first reference optimization strategy considered in this work is based on the MB family of PMFs, which as mentioned are known to be optimal in the linear region of transmission. The MB PMF is defined as $P_X(X) \propto \exp(-\lambda|X|^2)$, where $\lambda$ is a scale factor. It has been pointed out that the effective received SNR $SNR_{eff}$ in the nonlinear region is dependent on the modulation format through the 4th and 6th order moments of the constellation \cite{Julian}. In particular, constellations with high peak-to-average power ratio (PAPR) result in reduced $SNR_{eff}$ and overall degraded performance. To further exemplify this fact, two sets of MB PMFs are studied. One is optimal for an AWGN channel with an SNR, equivalent to the $SNR_{eff}$ with uniform input distribution to the channel at the respective launch power. Such constellations typically exhibit high PAPR. The other is optimized in a brute-force manner search of the parameter $\lambda$ at each launch power. In this case, $\lambda$ is allowed to be negative, i.e., resulting in high-energy points having higher probability than low-energy points, thus reducing the PAPR.

The third reference strategy is based on an MD ball constellation, known to be beneficial in the nonlinear region of transmission for constellations of high dimensionality \cite{Golani, PatrickMDball, Dar} and with a relatively large number of points. Similar to \cite{PatrickMDball}, the MD ball is constructed by pruning a 4D 16PAM constellation (256QAM per polarization) to the subset of $n_{ball}$ points of lowest energy, resulting in a low PAPR.   

\section{Results}
The AIRs are given in Fig.~1(b). For this distance, $64^2$QAM uniform outperforms $16^2$QAM uniform, and is thus used as reference for the rest of the results. When the MB is selected based on the received $SNR_{eff}$ (red circles), shaping gain is achieved in the linear region of transmission as expected. However, performance is highly penalized when the power is increased, resulting in a loss already before the optimum launch power. The loss is recovered by properly optimizing the MB (black triangles), resulting in a MB with $\lambda=0$ at the optimal launch power (equivalent to uniform distribution), and negative $\lambda$ at high powers. For the latter, high probability is assigned to high-energy points, resulting in the amplitudes being scaled down (for the same average power), and thus lower PAPR. 

Optimization with the proposed algorithm allows to increase the optimal launch power by 1 dB and achieve $\approx0.25$ and $\approx0.13$ bits/4D of gain when optimizing $64^2$QAM and $16^2$QAM constellations, respectively. Similar gain is achieved with the MD ball \cite{PatrickMDball} up to the optimal launch power. However, the latter employs a much larger constellation of 8192 points. The proposed algorithm effectively prunes points by assigning probability 0 to them, for a resulting optimum constellation with a total of 784 points, spread into 3 distinct amplitudes (see Fig.~1(c)). For high powers, the proposed algorithm converges to a single amplitude constellation, which gains $\approx0.5$ bits/4D w.r.t. the uniform, and similar w.r.t. the MD ball, which is now highly penalized. Finally, in Fig.~1(d), the gain w.r.t. $64^2$QAM uniform is given for the proposed algorithm, and for the MD ball as a function of the number of points allowed in the constellation. We see that when the constellation is constrained (e.g. due to the complexity, associated with detection and decoding), the MD ball is penalized even further, and at small constellations is even outperformed by uniform distribution. It is noted that the proposed algorithm does not explicitly search for small sized constellation, but naturally converges to a few-amplitude, few-points one. 
\begin{figure}[!t]

  \begin{minipage}{.7\textwidth}
  \subfigure{
  \includegraphics[width = 8.4cm ]{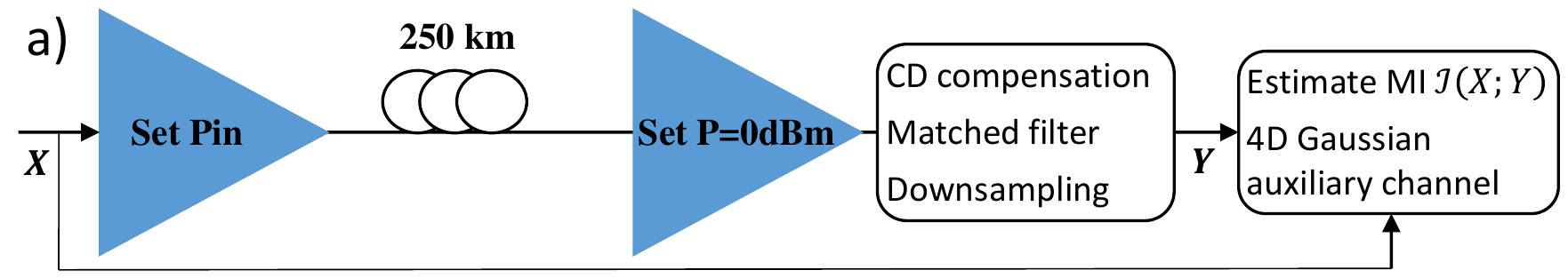}
  }
  \\
  \subfigure{
  \includegraphics[width = 4.4cm ]{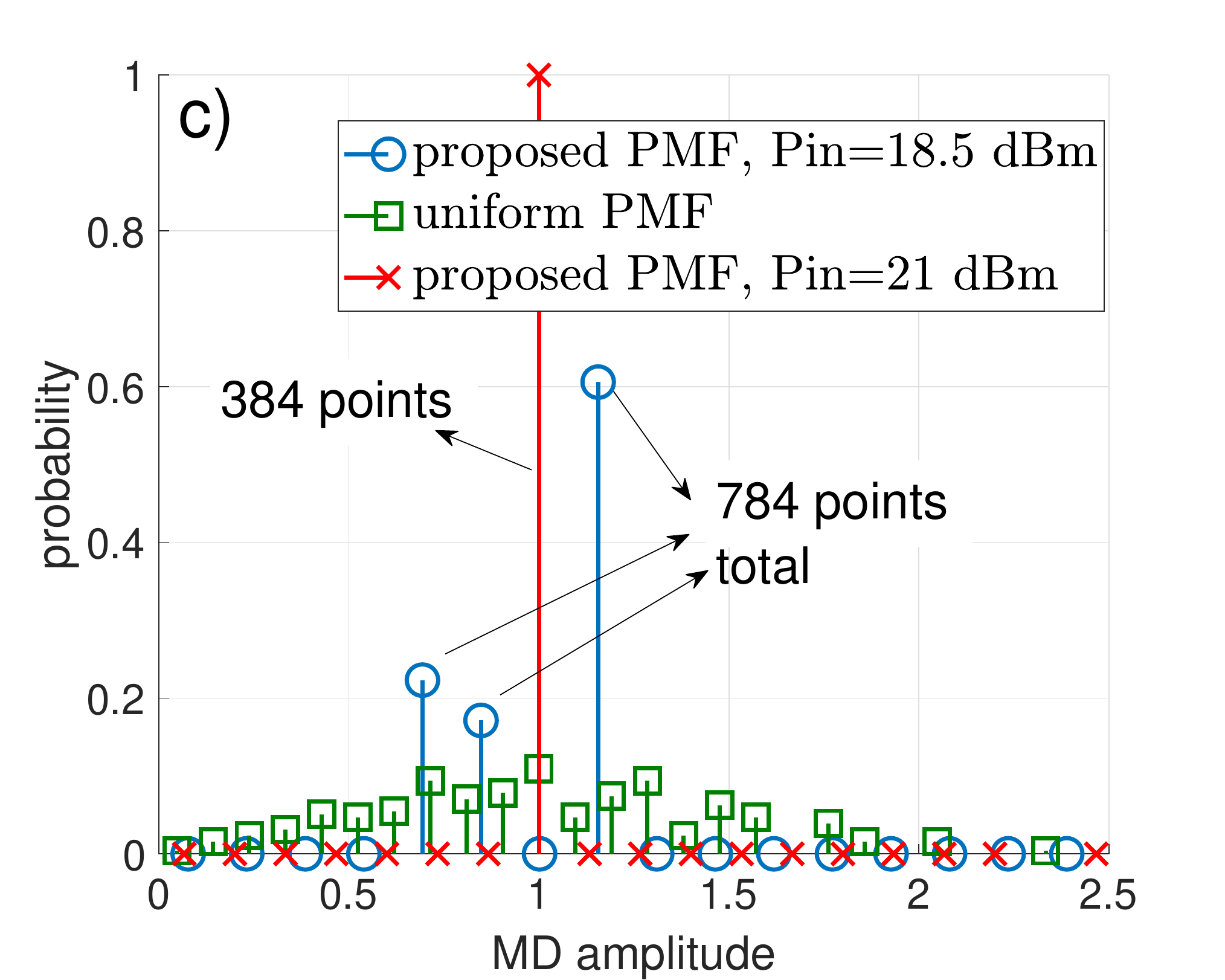}
  \includegraphics[width = 4.4cm ]{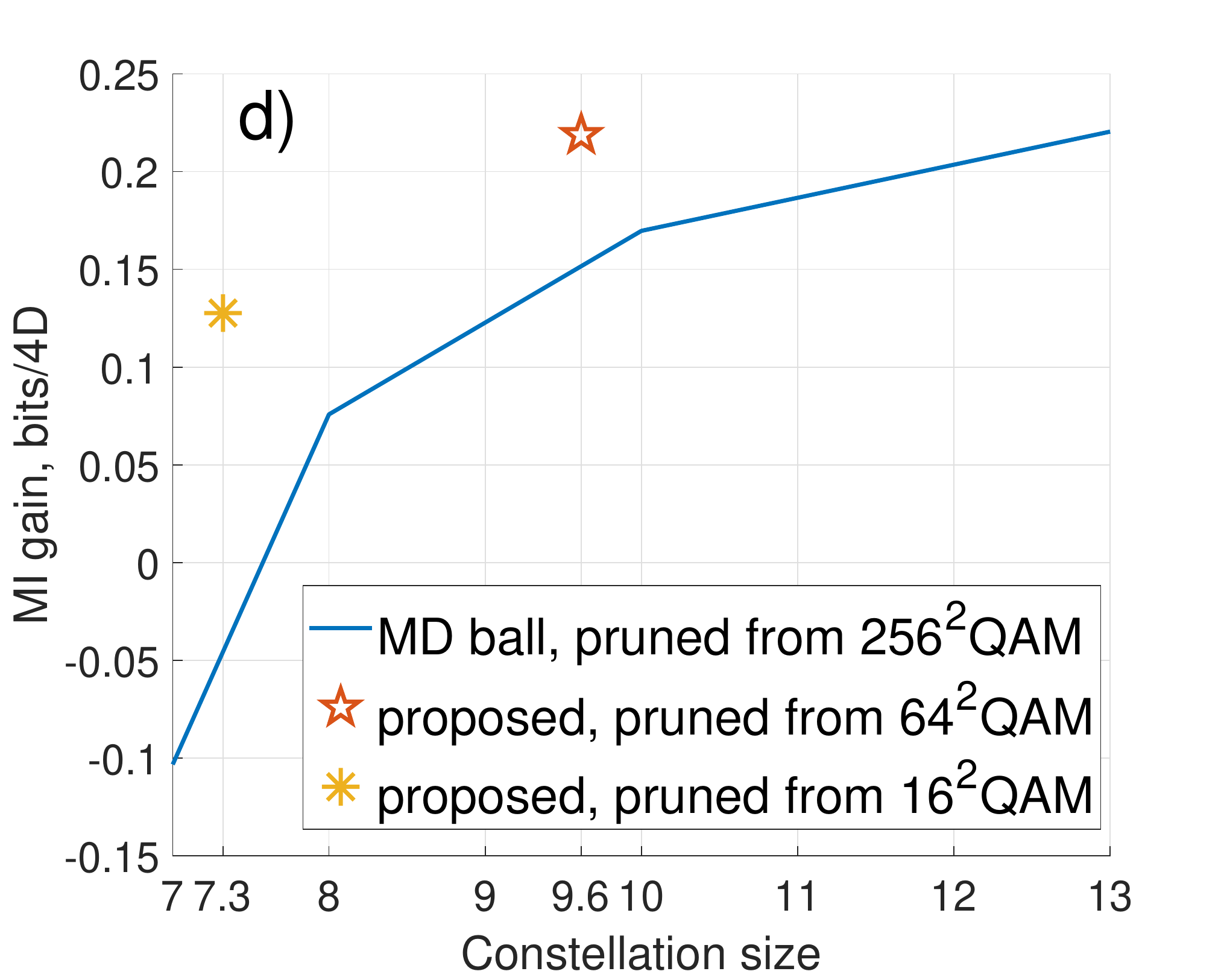}
  }
  \end{minipage}
  \begin{minipage}{.3\textwidth}
  \subfigure{
  \hspace{-2.6cm}
  \includegraphics[width = 7.3cm ]{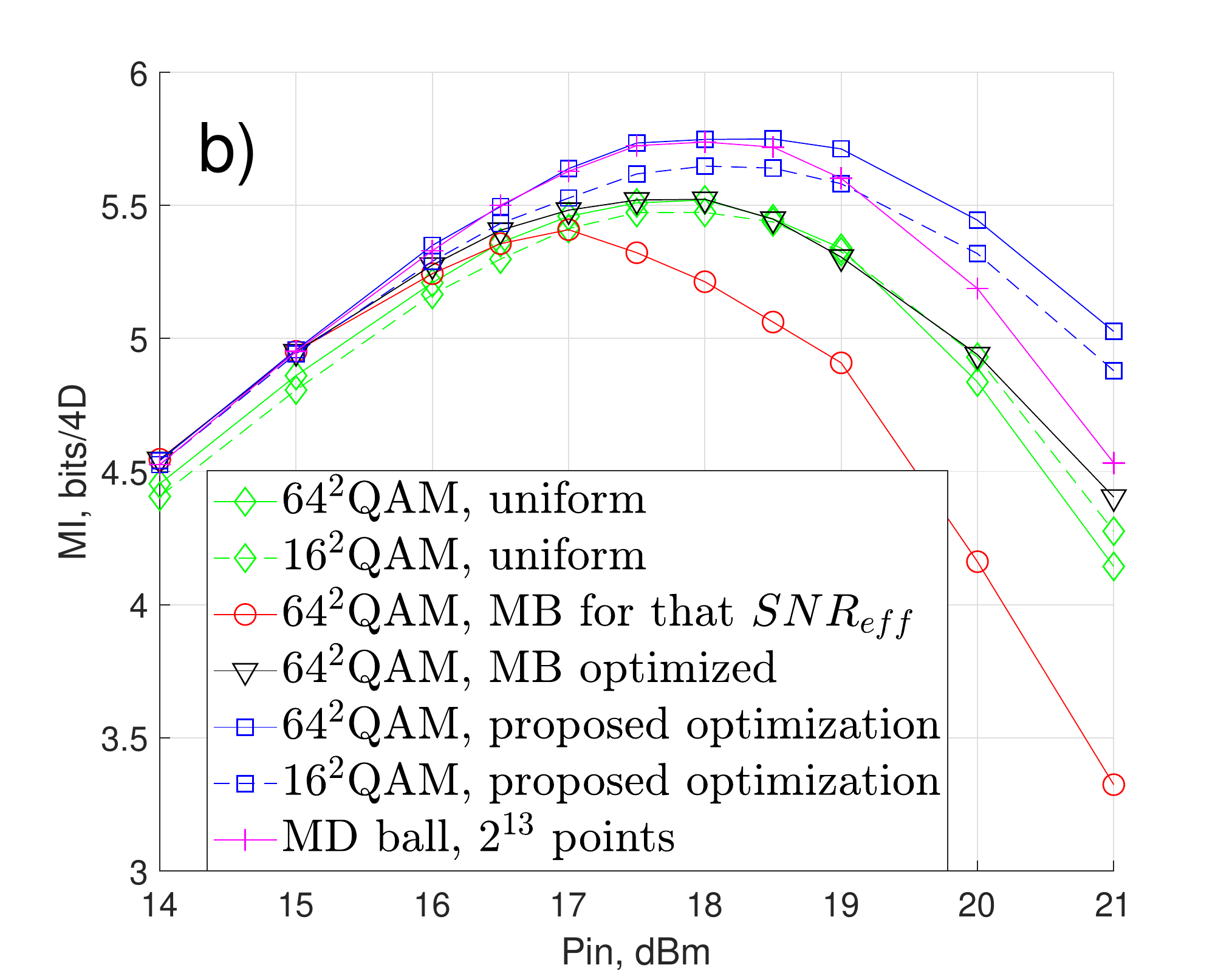}
  }
  \end{minipage}
  \caption{Results of the proposed optimization. \textbf{a):} Simulation setup; \textbf{b):} MI vs launch power for the studied optimization strategies; \textbf{c):} amplitude PMFs for the studied optimization strategies at the optimal launch power ($P_{in}=18.5$ dBm) and in the highly nonlinear region ($P_{in}=21$ dBm). In the first case, the proposed algorithm assigns 0 probability to most amplitudes, leaving a 3-amplitude constellation with a total of 784 points (pruned from 4096). In the second case, the proposed algorithm converges to a single-amplitude constellation; \textbf{d):} Achieved gain as a function of the number of points in the constellation. MD ball is heavily penalized when the constellation is small, resulting in a net loss w.r.t. uniform PMF.}
  \vspace{-.5cm}
  \label{fig:results}
  \end{figure}

\section{Conclusion}
An algorithm was proposed to optimize the PMF of multi-dimensional constellations. It was shown that in the nonlinear region, optimized PMFs exhibit a few-amplitude shape (single amplitude in the very high powers), which results in relatively small optimized constellation size. Both the lowest and highest amplitudes in the PMF are suppressed. This is in contrast with previous assumptions on the optimality of the PMF, targeting to use the points of lowest energy. In particular when the constellation is constrained to a small size, the proposed algorithm significantly outperforms the latter strategy. Employing an explicit constellation size constraint in the optimization can be expected to be further beneficial in keeping the same gains, while reducing the constellation size, thus keeping the complexity of the receiver low. 

\section*{Acknowledgments}
This work was supported by the DNRF Research Centre of Excellence, SPOC, ref. DNRF123.

\end{document}